\documentclass[11pt,leqno]{article}
\usepackage{amsfonts}
\usepackage{amsmath, amsthm, amsfonts, amssymb, color}
\usepackage{mathrsfs}
\usepackage{color}
\theoremstyle{definition}

\newcommand{\scr}[1]{\mathscr #1}
\definecolor{wco}{rgb}{0.5,0.2,0.3}

\numberwithin{equation}{section} \theoremstyle{remark}

\newcommand{\ua}{\uparrow}

\title{{\bf  $\Phi$-Entropy Inequality and Invariant Probability Measure for SDEs with Jump}\footnote{ Supported in
 part by  NNSFC(11131003), SRFDP, the Fundamental Research Funds for the Central Universities.} }
\author{
{\bf     Feng-Yu Wang  }\\
\footnotesize{ School of Mathematical Sciences,
Beijing Normal
University, Beijing 100875, China}\\
 \footnotesize{ Department of Mathematics,
Swansea University, Singleton Park, SA2 8PP, United Kingdom}\\
\footnotesize{  wangfy@bnu.edu.cn, F.-Y.Wang@swansea.ac.uk}}
\begin{document}
\allowdisplaybreaks
\def\R{\mathbb R}  \def\ff{\frac} \def\ss{\sqrt} \def\B{\mathbf
B}
\def\N{\mathbb N} \def\kk{\kappa} \def\m{{\bf m}}
\def\ee{\varepsilon}
\def\dd{\delta} \def\DD{\Delta} \def\vv{\varepsilon} \def\rr{\rho}
\def\<{\langle} \def\>{\rangle} \def\GG{\Gamma} \def\ggg{\gamma}
  \def\nn{\nabla} \def\pp{\partial} \def\E{\mathbb E}
\def\d{\text{\rm{d}}} \def\bb{\beta} \def\aa{\alpha} \def\D{\scr D}
  \def\si{\sigma} \def\ess{\text{\rm{ess}}}
\def\beg{\begin} \def\beq{\begin{equation}}  \def\F{\scr F}
\def\Ric{\text{\rm{Ric}}} \def\Hess{\text{\rm{Hess}}}
\def\e{\text{\rm{e}}} \def\ua{\underline a} \def\OO{\Omega}  \def\oo{\omega}
 \def\tt{\tilde} \def\Ric{\text{\rm{Ric}}}
\def\cut{\text{\rm{cut}}} \def\P{\mathbb P} \def\ifn{I_n(f^{\bigotimes n})}
\def\C{\scr C}      \def\aaa{\mathbf{r}}     \def\r{r}
\def\gap{\text{\rm{gap}}} \def\prr{\pi_{{\bf m},\varrho}}  \def\r{\mathbf r}
\def\Z{\mathbb Z} \def\vrr{\varrho} \def\ll{\lambda}
\def\L{\scr L}\def\Tt{\tt} \def\TT{\tt}\def\II{\mathbb I}
\def\i{{\rm in}}\def\Sect{{\rm Sect}}  \def\H{\mathbb H}
\def\M{\scr M}\def\Q{\mathbb Q} \def\texto{\text{o}} \def\LL{\Lambda}
\def\Rank{{\rm Rank}} \def\B{\scr B} \def\i{{\rm i}} \def\HR{\hat{\R}^d}
\def\to{\rightarrow}\def\l{\ell}\def\iint{\int}
\def\EE{\scr E}
\def\A{\scr A}\def\bGG{{\bf\Gamma}}
\def\BB{\scr B}\def\Ent{{\rm Ent}}

\maketitle

\begin{abstract}  By using the $\Phi$-entropy inequality derived in   \cite{Wu, Ch} for  Poisson measures, the same type of inequality is established for   a class of stochastic differential
equations driven by purely jump L\'evy processes. The semigroup $\Phi$-entropy inequality for SDEs driven by Poisson point processes as well as a sharp result on the existence of invariant probability measures are also presented.
\end{abstract} \noindent
 \vskip 2cm

\section{Introduction}

Let $\Phi\in C([0,\infty))\cap C^2((0,\infty))$ be convex such that $\Phi(0)=0$ and the function
$$\Psi_\Phi(u,v):= \Phi(u)-\Phi(v)-\Phi'(v)(u-v),\ \ u,v\ge 0$$ is non-negative and convex. Typical examples of $\Phi$ include $\Phi(u)=u\log u$ and $\Phi(u)=u^p$ for $p\in [1,2].$

Let $(\EE,\D(\EE))$ be a Dirichlet form on $L^2(\mu)$ for a probability measure $\mu$. The $\Phi$-entropy inequality considered in \cite{Ch} is of type
\beq\label{PE}\Ent_\mu^\Phi(f):= \mu(\Phi(f))-\Phi(\mu(f))\le C\EE(\Phi'(f),f),\ \ f,\Phi'(f)\in\D(\EE),f\ge 0 \end{equation} for some constant $C>0$. This inequality is equivalent to (see \cite[Corollary 1.1]{Ch})
\beq\label{PD} \Ent_\mu^\Phi(P_t f)\le \e^{-t/C}\Ent_\mu^\Phi(f),\ \ t\ge 0, f\in \B_b^+,\end{equation} where $P_t$ is the associated Markov semigroup and $\B_b^+$ is the set of all bounded positive elements in $L^2(\mu)$. When $\Phi(u)= u\log u$, the inequality \eqref{PE} reduces to the modified log-Sobolev inequality studied in \cite{Wu, SM}.

In this paper, we investigate the $\Phi$-entropy inequality for
  the following stochastic differential equation (SDE) on $\R^d$:
\beq\label{E1.1} \d X_t =b(X_t)\d t +\si\d L_t,\end{equation} where $b: \R^d\to\R^d$ is $C^1$-smooth with bounded $\nn b$, $\si$ is an invertible $d\times d$-matrix, and $L_t$ is a purely jump L\'evy process on $\R^d$ with L\'evy measure $\nu$, i.e. $L_t$ is generated by
\beq\label{LL}\scr L_0 f:=  \int_{\R^d}\big[f(\cdot+z)-f-\<\nn f, z\>1_{\{|z|\le 1\}}\big]\nu(\d z),\ \ f\in C_b^2(\R^d). \end{equation} Since $b$ is Lipschitz continuous, for any initial data $x\in\R^d$  the equation \eqref{E1.1} has a unique solution $X_t(x)$ for   $t\in [0,\infty)$.   Let $P_t$ be the associated Markov semigroup, i.e.
$$P_t f(x):= \E f(X_t(x)),\ \ t\ge 0, f\in \B_b(\R^d), x\in\R^d,$$ where $\B_b(\R^d)$ is the set of all bounded measurable functions on $\R^d$.

When $P_t$ has an invariant probability measure $\mu$, we consider the corresponding (possibly non-sectorial) form
\beq\label{**}\EE(f,g):=-\int_{\R^d} f\scr Lg\,\d\mu,\ \ f,g\in C_0^2(\R^d),\end{equation} where $\scr L$ is the generator of $P_t$, i.e.
\beq\label{L} \scr L f=  \<\nn f, b\> +\int_{\R^d} \big[f(\cdot+\si z)-f-\<\nn f,\si z\> 1_{\{|z|\le 1\}}\big]\nu(\d z).\end{equation}

Let
$$\GG_{\Phi,\nu}(f)(x)=\int_{\R^d}\Psi_{\Phi}\big(f(x+\si z), f(x)\big)\nu(\d z),\ \ f\in \B_b^+(\R^d),$$ where $\B_b^+(\R^d)$ is the set of all positive elements in $\B_b(\R^d).$
Let $ C_{c,+}^2(\R^d)$ be the set of any $C^2$ positive function  on $\R^d$ which is constant outside a compact set. Then for any $f\in C_{c,+}^2(\R^d)$ we have $\int_{\R^d}\scr L \Phi(f)\d\mu=0$, so that   \eqref{L}  yields
\beg{equation}\label{W3}\beg{split} &\EE(\Phi'(f),f) := -\int_{\R^d} \Phi'(f) \scr L f\,\d\mu \\
&=\int_{\R^d}\d\mu \int_{\R^d} \Psi_\Phi\big(f(\cdot+\si z),f\big)\nu(\d z)- \int_{\R^d}\scr L \Phi(f)\d\mu\\
&= \int_{\R^d} \GG_{\Phi,\nu}(f)\d\mu.\end{split}\end{equation}
Thus, for the present model, the $\Phi$-entropy inequality \eqref{PE} reduces to
\beq\label{PE'} \Ent_\mu^\Phi(f)\le C \int_{\R^d} \GG_{\Phi,\nu}(f)\d\mu,\ \ f\in \B_b^+(\R^d).\end{equation}

\beg{thm}\label{T1.1} Assume that $\ff {\kk_1}{|z|^{d+\aa}}\le \ff{\nu(\d z)}{\d z}\le \ff {\kk_2}{|z|^{d+\aa}}$ for some constants $\kk_1,\kk_2>0$ and $\aa\in (0,2).$ Let $\ll_1,\ll_2\in\R$ such that
\beq\label{A0} \ll_1|v|^2\le \<\si^{-1} (\nn b(x))\si v, v\>\le \ll_2 |v|^2, \ \ x,v\in\R^d.\end{equation}
\beg{enumerate}\item[$(1)$] For any $T>0$ and $f\in \B_b^+(\R^d)$,
$${\rm Ent}_{P_T}^\Phi(f):= P_T\Phi(f)- \Phi(P_Tf)\le \ff{\kk_2\big(\exp[\ll_2(d+\aa)T-\ll_1Td]-1\big)}{\kk_1\big(\ll_2(d+\aa)-\ll_1 d\big)}P_T \GG_{\Phi,\nu}(f).$$
\item[$(2)$] If $\ll_2(d+\aa)<\ll_1 d,$  then $P_t$ has a unique invariant probability measure $\mu$ and  $\eqref{PE'}$ holds for $C:= \ff {\kk_2} {\kk_1(\ll_1d-\ll_2(d+\aa))}.$
     \end{enumerate}\end{thm}

 The following result partly  extends Theorem \ref{T1.1}  to the case where the L\'evy process $L_t$  merely has  large (e.g. $\rr= 1_{[1,\infty)}$) or small (e.g. $\rr=1_{(0,1]}$) jumps. In particular, \eqref{*W} holds in the situation of Theorem \ref{T1.1}(2).  

\beg{thm}\label{T1.2} Let $\ff{\kk_1\rr(|z|)}{|z|^{d+\aa}}\le \ff{\nu(\d z)}{\d z}\le \ff{\kk_2\rr(|z|)}{|z|^{d+\aa}}$   for some constants $\kk_1,\kk_2>0$ and some non-negative measurable function $\rr$ on $(0,\infty).$ Assume that $\eqref{A0}$ holds. \beg{enumerate} \item[{\rm (I)}] If $\ll_2\le 0$ and $\rr$ is decreasing, then  assertions $(1)$ and $(2)$ in  Theorem $\ref{T1.1}$ hold. In particular, if $\ll_2(d+\aa)<\ll_1d $  then $P_t$ has a unique invariant probability measure $\mu$ such that
 \beq\label{*W} \Ent_\mu^\Phi(P_t f)\le \exp\Big[-\ff{\kk_1(\ll_1 d-\ll_2(d+\aa))}{\kk_2}t\Big]\Ent_\mu^\Phi(f),\ \ t\ge 0, f\in \B_b^+(\R^d).\end{equation} \item[{\rm (II)}] If $\ll_1\ge 0$ and $\rr$ is increasing, then the assertion $(1)$ in Theorem $\ref{T1.1}$ holds. \end{enumerate}\end{thm}

\paragraph{Remark 1.1.} (1) We would like to mention a nice entropy inequality   derived recently in
\cite{JW} for non-local Dirichlet forms.  Let $\mu(\d
x):=\e^{-V(x)}\d x$ be a probability measure on $\R^d$ and let $\rr$
be a positive function on $(0,\infty)$ such that
\beq\label{AB}c:=\inf_{x,y\in \R^d}
\rr(|x-y|)\{\e^{V(x)}+\e^{V(y)}\}>0,\end{equation} then
\beg{equation*}\beg{split} &\Ent_\mu(f):= \mu(f\log f)-\mu(f)\log\mu(f)\\
&\le \ff 1 c \int_{\R^d} \mu(\d x)\int_{\R^d}
\Big\{(f(x+z)-f(z))\log\ff{f(x+z)}{f(x)}\Big\}\rr(|z|)\d z,\ \
f\in\B_b^+(\R^d).\end{split}\end{equation*}  Since $\Psi_{\log} (u,v)\le (u-v)\log\ff uv$ for $u,v>0$, this inequality follows from the corresponding $\Phi$-entropy inequality with $\Phi(r)= r\log r.$  But, in general  this result is incomparable with ours for $\Phi(r)=r\log r$. In our case   the invariant probability measure
of $P_t$ (if exists) is not explicitly known, so that the condition (\ref{AB}) is hard to verify.   Moreover, condition
(\ref{AB}) implies that $\nu(\d z):= \rr(|z|)\d z$ has full support
on $\R^d$ which does not apply to the situations of Theorem \ref{T1.2} if $\rr$ is not strictly positive on $(0,\infty)$.

(2) When $b(x)=-x$ and $\nu(\d z)= N(z)\d z$ for $N\ge 0$ satisfying
$$\int_1^\infty N(sz)s^{d-s}\d s\le CN(z)$$ for some constant $C>0$, the $\Phi$-entropy inequality \eqref{PE'} was proved in \cite{GI}. This condition is satisfied for $\nu$ given in Theorem \ref{T1.1} but fails in the situation of Theorem \ref{T1.2}(1) for e.g. $\rr(s)=1_{\{s\le 1\}}.$ Moreover, as shown in (d) in the proof of Theorem \ref{T1.2}, to deduce the exponential convergence from \eqref{PE'} an approximation argument should be  included in the proof of \cite[Theorem 1]{GI}, since the formula in \cite[Lemma 1]{GI} only makes sense for functions $w_1,w_2$ such that $w_1 L w_2\in L^1(u_\infty \d x)$. In general, the form $(\EE, C_0^2(\R^d))$ given in \eqref{**} does not provide a Dirichlet form, so that the equivalence between
\eqref{PE} and \eqref{PD} for Dirichlet forms does not apply.

\

Next, partly for the proof of Theorem \ref{T1.1}(2), we consider the existence of invariant probability measures for the following more general SDE:
\beq\label{3.1} \d X_t= b(X_t) \d t +\si_1(X_t)\d W_t +\si_2(X_{t-})\d   L_t,\end{equation} where $b: \R^d\to\R^d$ and $\si_1,\si_2: \R^d\to\R^d\otimes\R^d$ are locally Lipschitz continuous, $  L_t$ is the L\'evy process in \eqref{E1.1},    and $W_t$ is a $d$-dimensional Brownian motion independent of $L_t$. Then \eqref{3.1} has a unique solution up to the life time.

Although   the existence of invariant probability measures for SDEs with jumps has been investigated in the literature, we did not find any existing result which  directly applies  to the framework in Theorem \ref{T1.1}. For instance, in \cite[Theorem 4.5]{ABW} it is assumed that $\int_{\R^d} |z|^2\nu(\d z)<\infty$, while in \cite{AW} the L\'evy process is assumed to be the $\aa$-stable process and $b(x)$ is a perturbation by $-\ggg x$ for some constant $\ggg>0,$ see also \cite{CR,CZ} for the study of semilinear SPDEs with jump.  We aim to present a new result which is sharp in terms of the L\'evy measure and, in particular,  implies the existence of invariant probability measure in the situation of  Theorem \ref{T1.1}(2).

\beg{thm}\label{T3.1} Let $B\in C^1([0,\infty))$ be strictly positive. For any $\vv>0$ let
$$\tt B_\vv(x)= \sup\Big\{\ff {B(r)-rB'(r)}{2B(r)^2(1+r)}:\ r\ge 0, |r-|x||\le \vv \|\si_2(x)\|\Big\},\ \ x\in\R^d.$$ If  there exists $\vv\in (0,1]$ such that either \beg{equation}\label{C1}\beg{split} A_\vv  :=  \limsup_{|x|\to\infty} \bigg\{&\ff{\<b(x),x\>+{\rm Tr}(\si_1\si_1^*)(x)}{B(|x|)(|x|+1)}+\ff{\|\si_2(x)\|\int_{\{\vv<|z|\le 1\}}|z|\nu(\d z)}{B(|x|)} \\
&-\ff{|\si_1^*(x)x|^2B'(|x|)}{B(|x|)^2(1+|x|)|x|}+\int_{\{|z|>\vv\}}\nu(\d z) \int_{|x|}^{|x|+\|\si_2(x)\|\cdot|z|}\ff{\d s}{B(s)}\\
 & + \|\si_2(x)\|^2\tt B_\vv(x)\int_{\{|z|\le \vv\}}|z|^2\nu(\d z)
 \bigg\}=-\infty,\end{split}\end{equation}   or
 \beq\label{C2} A_\vv<0\ \text{and}\ \int_0^\infty \ff{\d s}{B(s)}=\infty,\end{equation} then the
   solution to $\eqref{3.1}$ is non-explosive and the associated Markov semigroup has an invariant probability.  \end{thm}

The following is  a consequence of Theorem \ref{T3.1}, which provides some more explicit sufficient conditions for the existence of invariant probability measures.

\beg{cor} \label{C3.2} Assume that   for some  $\theta\in\R$
\beq\label{C} D:=\limsup_{|x|\to\infty}  \bigg\{\ff{\<b(x),x\>+ {\rm Tr}(\si_1\si_1^*)(x)}{  (1+|x|)^{1+\theta}}-\ff{\theta|\si_1^*(x)x|^2}{|x|(1+|x|)^{\theta+2}}\bigg\}<0,\end{equation} and that
$$\Theta:= \limsup_{|x|\to\infty} \ff{\|\si_2(x)\|}{|x|}<\infty.$$
Then the solution to $\eqref{3.1}$ is non-explosive and the associated Markov semigroup has an invariant probability measure in each of the following three situations:
\beg{enumerate}\item[$(1)$] $\theta>1$.   \item[$(2)$]   $\theta= 1, \int_{\{|z|\ge 1\}}\log(1+|z|)\nu(\d z)<\infty$,   and there exists $\vv\in (0,\Theta^{-1})\cap (0,1]$ such that
\beq\label{20} \ff{\vv^2\Theta^2\int_{\{|z|\le \vv\}}|z|^2\nu(\d z)}{2(1-\vv\Theta)^2}+\Theta\int_{\{\vv<|z|\le 1\}}|z|\nu(\d z) + \int_{\{|z|>\vv\}} \log (1+\vv\Theta |z|)\nu(\d z) <-D.\end{equation}
\item[$(3)$] $\theta\in (0,1), \|\si_2\|$ is bounded, and $\int_{\{|z|\ge 1\}} |z|^{1-\theta}\nu(\d z) <\infty.$
\item[$(4)$] $\theta\in (-\infty,1),  \int_{\{|z|\ge 1\}} |z|^{1+\theta^-} \nu(\d z) <\infty,$ and
\beq\label{21} \limsup_{|x|\to\infty}  \ff{\|\si_2(x)\|}{|x|^\theta} \int_{\{|z|> 1\}} |z|  \,\nu(\d z)<-D.\end{equation} \end{enumerate} \end{cor}

 Note that when \eqref{C} holds with $\theta=1$ and $\lim_{|x|\to\infty} \ff{\|\si_2(x)\|}{|x|}=0$, Corollary \ref{C3.2} implies the existence of the invariant probability measure provided
 \beq\label{SH}\int_{\{|z|\ge 1\}}\log (1+|z|)\nu(\d z)<\infty.\end{equation} According to \cite[Theorems 17.5 and 17.11]{SA}, \eqref{SH} is   sharp (i.e. sufficient and necessary) for the purely jump Ornstein-Uhlenbeck process (i.e. $\si_1=0, \si_2=I,  b(x)=-x$) to have invariant probability measure.
When $\theta\in (0,1), \si_1=0, \si_2=I$ and $b(x)= -x |x|^{\theta-1},$ we would believe that the condition $\int_{\{|z|\ge 1\}} |z|^{1-\theta}\nu(\d z) <\infty$ in case (3) is also sharp for the existence of the invariant probability measure. However, in this case the distribution of the solution is no longer infinitely divisible, so that the proof of  \cite[Theorem 17.11]{SA} does not apply.

 \

The remainder of the paper is organized as follows. In Section 2, by using the $\Phi$-entropy inequality derived in  \cite{Wu} and \cite{Ch} for  Poisson measures,   we prove a   result on the semigroup $\Phi$-entropy inequality for SDEs driven by Poisson point processes. In Section 3 we
prove Theorem \ref{T3.1} and Corollary \ref{C3.2}.   Finally, proofs of Theorems \ref{T1.1} and \ref{T1.2} are presented in Section 4.

\section{The semigroup $\Phi$-entropy inequality}

 Let   $N(\d t,\d z)$ be  a Poisson point process on $\R^d$ with compensator $\d t\,\nu(\d z)$, where $\nu$ is a $\si$-finite measure on $\R^d$.
Then for any $T>0$, $1_{[0,T]}(t)N(\d t,\ dz)$ is a random variable on the configuration space
$$\bGG_T:= \Big\{\ggg=\sum_{i=1}^n \dd_{(s_i,z_i)}: n\in \Z_+\cup\{\infty\}, (s_i,z_i)\in [0,T]\times\R^d\Big\} $$     equipped with the $\si$-field
induced by $\{\ggg\mapsto\ggg(A): A\in \B([0,T]\times\R^d)\},$ where  $\B([0,T]\times\R^d)$ is the Borel $\si$-field on $[0,T]\times\R^d$ and $\dd_{(s_i,x_i)}$ stands
for the Dirac measure at point $(s_i,x_i)$. The distribution of $1_{[0,T]}(t)N(\d t,\ dz)$ is the Poisson measure with intensity $ \d t\,\nu(\d z)$
on $[0,T]\times\R^d$.

 Let
$$a: [0,\infty)\times \R^d\to\R^d$$ be measurable such that  for every $s\ge 0$, $a_s$ is invertible and
\beq\label{*C1}  \int_{[0,t]\times\R^d} (1\land |a_s(z)|^2)\d s\nu(\d z)<\infty,\ \ t\ge 0.\end{equation}
Let
$$ \check N_a(\d t,\d z)= N(\d t,\d z)- 1_{\{|a_t(z)|\le 1\}}\d t \nu(\d z).$$
Then the stochastic integral
$$\int_{[0,t]\times \R^d}  a_s(z) 1_{\{|a_s(z)|\le 1\}}\check N(\d s,\d z),\ \ t\ge 0$$ is well defined (see e.g. \cite[page 36-37]{ST}). Moreover, since (\ref{*C1}) implies that  $\P$-a.s.,
$$1_{\{s\in [0,t], |a_s(z)|>1\}}N(\d s,\d z)\in \bGG_t^0:=\big\{\ggg\in \bGG_t: \ggg([0,t]\times\R^d)<\infty\big\},$$    the stochastic integral
$$\int_{[0,t]\times \R^d}  a_s(z) \check N_a(\d s,\d z),\ \ t\ge 0$$ is well defined as well.

Now, consider the
 following equation on $\R^d$:
 \beq\label{E1}   \d X_t=  b_t(X_t)\d t+\int_{\R^d} a_t(z)\check N_a(\d t,\d z),\ \ t\ge 0,\end{equation} where   $b:
 [0,\infty)\times\R^d\to\R^d$ is measurable such that   $b_t$ is
 Lipschitz continuous for every $t\ge 0$ and the Lipschitz constant is locally bounded
 in $t$. It is standard that for any $x\in\R^d$, this equation has a unique solution $X_t(x)$ with $X_0=x$, see e.g. \cite[Theorem 17]{ST}.

 Let
 $$P_t f(x)= \E f(X_t(x)),\ \ t\ge 0, x\in\R^d, f\in \B_b(\R^d).$$ We
 aim to establish the   $\Phi$-entropy inequality for $P_T$.
 To state our main result, we introduce the following equation driven by $\check N_a+\dd_{(s,z)}$ for $(s,z)\in (0,\infty)\times \R^d$:
\beq\label{XS}  X_t^{s,x}(z)= x+ \int_0^t b_r(X^{s,x}_r(z))\d r + \int_{[0,t]\times \R^d} a_r(y) \{\check N_a+\dd_{(s,z)}\}(\d r,\d y),\ \ t\ge 0.\end{equation}

\beg{thm}\label{T2.1}   For fixed  $T>0$ and $x\in\R^d$, let
$$\psi_{s}(z)= a_s^{-1}\big(X_T^{s,x}(z)-X_T(x)\big),\ \ s\in (0,T], z\in\R^d.$$
If $\nu\circ\psi_{s}^{-1}$ is absolutely continuous w.r.t. $\nu$
such that \beq\label{CC}\xi_s:= {\rm ess}_{\P\times\nu}
\ff{\d\nu\circ\psi_{s}^{-1}}{\d\nu}<\infty,\ \ s\in
(0,T],\end{equation}then
$$\Ent_{P_T}^\Phi(f)(x)\le \E\int_0^T \xi_t\d t \int_{\R^d} \Psi_\Phi\big(f(X_T(x)+a_t(z)), f(X_T(x))\big) \nu(\d z),\ \ x\in\R^d, f\in\B_b^+(\R^d).$$  \end{thm}

Throughout this section, we fix $T>0$ and $x\in\R^d,$ and simply denote
  $$N_T:=1_{[0,T]}(t)N(\d s,\d z).$$  To prove Theorem
\ref{T2.1}, we shall use  the following $\Phi$-entropy inequality for the Poisson point process on $[0,T]\times\R^d$:
 \beq\label{Wu2}  \beg{split} &\E(\Phi\circ F (N_T)) -\Phi(\E
F(N_T))\\
&\le \E\int_{[0,t]\times\R^d}\Psi_\Phi\big(F(N_T+\dd_{(s,z)}),
F(N_T)\big) \d s\nu(\d z),\ \ F\in
\B_b^+(\bGG_T).\end{split}\end{equation} This inequality was first proved by Wu \cite{Wu} for $\Phi(u)= u\log u$, and as explained in  \cite[\S 5.1]{Ch} that Wu's proof also applies to general $\Phi$ considered in the paper.

According to the
inequality \eqref{Wu2}, to prove Theorem \ref{T2.1}   we need   to formulate
$X_T(x)+a_t(z)$ using $N_T+ \dd_{(\tau,\xi)}$ for some $\xi\in\R^d$ and
$\tau\in [0,T].$ To this end, we let $F: [0,T] \times \bGG_T\to\R^d$ be
measurable such that $X_t(x) = F_t(N_T), t\in [0,T]$.
Then we  would  suggest
that $Y_t:= F_t(N_T+\dd_{(\tau,\xi)})$   solves the equation
$$  Y_t=x+ \int_0^t b_s(Y_s)\d s +   \int_{[0,t]\times \R^d} a_s(z)\{\check N_a+\dd_{(\tau,\xi)}\}(\d s,\d z),\ \ t\in [0,T].$$ Thus, taking $\xi$ and $\tau$ such that $a_t(z)= Y_T-X_T(x)$, we obtain
$$X_T(x)+a_t(z)= F_T(N_T+\dd_{(\tau,\xi)}).$$ However, since
$X_\cdot(x)=F_\cdot(N_T)$ holds on $[0,T]$ merely $\P$-a.s.,   to
make this argument rigorous we need  to verify the quasi-invariance for the transform $N_T\to N_T+\dd_{(\tau,\xi)}$,
which is ensured by the following  Girsanov type theorem, see   \cite{W} for a similar result for L\'evy processes.

\beg{lem}\label{L2.1} Let $g$ be a strictly positive function on
$[0,T]\times\R^d$ such that $\nu^g(\d s,\d z):= g(s,z)\d s\nu(\d z)$ is a probability measure on $[0,T]\times\R^d.$ Let
$$N_T(g)=\int_{[0,T]\times\R^d} g(s,z)N(\d s,\d z).$$ Moreover, let
$(\tau,\xi)$ be a random variable independent of $N_T$ and with
distribution $\nu^g$. Then
$$R:=\ff 1{g(\tau,\xi)+N_T(g)}$$ is a strictly positive probability
density w.r.t. $\P$ such that  the distribution of
$N_T+\dd_{(\tau,\xi)}$ under $\d\Q:=R\d\P$ coincides with that of
$N_T$ under $\P$. \end{lem}

\beg{proof} Let $\pi$ be the Poisson measure with intensity   $\d s\nu(\d z)$ on $[0,T]\times\R^d$. Then $\pi\times\nu^g$ is the
distribution of $(N_T,\tau,\xi)$. By the Mecke formula for the Poisson measure
(see (3.1) in \cite{M}), for any $F\in \B_b^+({\bf \ggg}_T)$ we have
\beg{equation*}\beg{split} &\E\big\{RF(N_T+\dd_{(\tau,\xi)})\} =
\int_{{\bf \ggg}\times[0,T]\times\R^d}\ff{F(\ggg+\dd_{(s,z)})
g(s,z)}{(\ggg+\dd_{(s,z)})(g)}\,\pi(\d\ggg)\d s\nu(\d z)\\
&= \int_{\bGG}
\ff{F(\ggg)\ggg(g)}{\ggg(g)}\,\pi(\d\ggg)=\pi(F).\end{split}\end{equation*}
Therefore, $\Q:= R\d\P$ is a probability measure, and the
distribution of $N_T+\dd_{(\tau,\xi)}$ under $\Q$ coincides with
that of $N_T$ under $\P$.
\end{proof}

\beg{proof}[Proof of Theorem \ref{T2.1}]   Let $F:   \bGG_T\to\R^d$ be measurable such
that $X_T(x)= F(N_T)$. We intend to prove
\beq\label{*}X_T^{s,x}(z)= F(N_T+\dd_{(s,z)}),\ \ \P\times\d s\times \nu(\d z)
\text{-a.e.}\end{equation} To this end, for $ g\in \B^+([0,T]\times\R^d)$   in Lemma
\ref{L2.1}, consider the product probability space:  $$\bar\OO=
\OO\times[0,T]\times\R^d,\ \ \bar\P(\d\oo,\d s,\d z)=g(s,z)
\P(\d\oo)\d s\nu(\d z).$$ Let $\bar N= (N,\tau,\xi)$ be defined by
$$  N(\oo,s,z)=  N(\oo),  \ \ \xi(\oo,s,z)=z,\ \ \tau(\oo,s,z)=s,\ \ (\oo,s,z)\in\bar\OO.$$
Then under $\bar\P$  the random variable $(\tau,\xi)$ is independent
of $N_T$ and has distribution $\nu^g(\d s,\d z):=g(s,z)\d s\nu(\d z) $.  Let $R$ be in Lemma \ref{L2.1}. Then the distribution of
$N_T+\dd_{(\tau,\xi)}$ under $\Q$  coincides with that of $N_T$
under $\bar \P$ (equivalently, under $\P$).  Thus, by the weak uniqueness
of solutions to (\ref{E1}),
  the
distribution of $(N_T+\dd_{(\tau,\xi)},Y_T)$ under $\Q$ coincides
with that of $(N_T,X_T(x))$ under $\P$. In particular,  the
distribution of $Y_T-F(N_T+\dd_{(\tau,\xi)})$ under $\Q$ coincides
with that of $X_T(x)- F(N_T)$ under $ \P$. Since $X_T(x)=F(N_T)\ \P$-a.s.,
this implies that
$$Y_T= F(N_T+\dd_{(\tau,\xi)}),\ \ \Q\text{-}a.s.$$ As $\Q$
is equivalent to $\bar\P$, it also holds $\bar\P$-a.s. Then (\ref{*}) follows   by noting that $Y_T= X_T^{\tau,x}(\xi)$ and $g>0$ such that $\bar\P$ is equivalent to $\P\times\d s\times\nu(\d z).$

Now, by (\ref{Wu2}) and (\ref{*}), for any $f\in
\B_b^+(\R^d)$ we have \beq\label{*2}\beg{split} \Ent_{P_T}^\Phi (f) &\le
\E\int_{[0,T]\times\R^d} \Psi_\Phi\big(f\circ
F(N_T+\dd_{(s,z)}), f\circ F(N_T)\big)  \d s\nu(\d z) \\
&= \E \int_{[0,T]\times\R^d} \Psi_\Phi\big(f(X_T^{s,x}(z)),
f(X_T(x))\big)\d s\nu(\d z).\end{split}\end{equation} Noting that
$X_T^{s,x}(z)= X_T(x) +a_s\circ\psi_{s}(z)$, it follows from (\ref{CC}) that
\beg{equation*}\beg{split} &\E\int_{\R^d} \Psi_\Phi\big(f(X_T^{s,x}(z)),
f(X_T(x))\big) \nu(\d z)\\
  &= \E\int_{\R^d} \Psi_\Phi\big(f(X_T(x)+a_s(z)),
f(X_T(x))\big) (\nu\circ\psi_{s}^{-1})(\d z)\\
&\le \E\int_0^T\xi_s\d s \int_{\R^d} \Psi_\Phi\big(f(X_T(x)+a_s(z)), f(X_T(x))\big)\nu(\d
z),\ \ s\in (0,T].\end{split}\end{equation*} Combining this with (\ref{*2}) we
finish the proof.
\end{proof}

\section{Proofs of Theorem \ref{T3.1} and Corollary \ref{C3.2}}

\beg{proof}[Proof of Theorem \ref{T3.1}]  Take $W(x)=\varphi(|x|)$, where   $\varphi(r):= \int_0^r \ff{s}{(1+s)B(s)}\d s,\ r\ge 0.$ Then $W\in C^2(\R^d)$. Let $\scr L$ be the generator of the solution $X_t$ to \eqref{3.1}. By the It\^o formula we have
\beg{equation}\label{AB}\beg{split} \scr L W(x)= &\<b(x), \nn W(x)\> + {\rm Tr}\big[(\si_1\si_1^*  \nn^2W)(x)\big] \\ &+\int_{\R^d}\big[W(x+\si_2(x)z)-W(x)-\<\nn W(x), \si_2(x) z\>1_{\{|z|\le 1\}}\big]\nu(\d z)\end{split}\end{equation} if the integral in the right hand side exists.  We observe that it suffices to prove that $\scr LW$ is a well defined locally bounded function with
\beq\label{C3} \beg{split} \scr LW(x)\le &\ff{\<b(x),x\>+{\rm Tr}(\si_1\si_1^*)(x)}{B(|x|)(|x|+1)}+\ff{\|\si_2(x)\|\int_{\{\vv<|z|\le 1\}}|z|\nu(\d z)}{B(|x|)} \\
&-\ff{|\si_1^*(x)x|^2B'(|x|)}{B(|x|)^2(1+|x|)|x|}+\int_{\{|z|>\vv\}}\nu(\d z) \int_{|x|}^{|x|+\|\si_2(x)\|\cdot|z|}\ff{\d s}{B(s)}\\
 & + \|\si_2(x)\|^2\tt B_\vv(x)\int_{\{|z|\le \vv\}}|z|^2\nu(\d z)  .\end{split}\end{equation}
In fact, by this and $A_\vv=-\infty$ we see that $-\scr LW$ is a compact function (i.e. $\{-\scr LW\le r\}$ is relatively compact for $r>0$). Therefore, by the It\^o formula we see that the solution is non-explosive with
$$\limsup_{t\to\infty} \ff 1 t \int_0^t \E(-\scr LW)(X_s)\d s <\infty,$$ which implies the existence of the invariant probability measure by a standard tightness argument. Moreover, if $\int_0^\infty \ff{\d s}{B(s)}=\infty$ and $A<0$, then $W$ is a compact function and by the It\^o formula the solution is non-explosive with $\E W(X_t)<\infty, t\ge 0.$ Thus, according to \cite[Theorem 4.1]{CD}, $A<0$   also implies that   the associated Markov semigroup has an invariant probability measure.  Below we prove that $\scr LW$ is locally bounded such that \eqref{C3} holds.
\beg{enumerate}
\item[(a)] It is easy to see that
\beg{equation*}\beg{split} &\big\{\<b, \nn W\> + {\rm Tr} (\si_1\si_1^*\nn^2W )\big\}(x)-\int_{\{\vv<|z|\le 1\}}\<\nn W(x), \si_2(x)z\>\nu(\d z)\\
&\le  \varphi'(|x|)  \Big(\ff{\<b(x),x\>+{\rm Tr}(\si_1\si_1^*)(x)}{|x|}-\ff{|\si_1^*(x)x|^2}{|x|^3}+\|\si_2(x)\|\int_{\{\vv<|z|\le 1\}}|z|\nu(\d z)  \Big)\\
&\quad + \varphi''(|x|)\ff{|\si_1^*(x) x|^2}{|x|^2}\\
&\le \ff{\<b(x),x\>+{\rm Tr}(\si_1\si_1^*)(x)}{B(|x|)(1+|x|)}+\ff{\|\si_2(x)\|\int_{\{\vv<|z|\le 1\}}|z|\nu(\d z)}{B(|x|)}\\
& \quad -\ff{|\si_1^*(x)x|^2B'(|x|)}{B(|x|)^2(1+|x|)|x|},\ \ x\in\R^d.\end{split}\end{equation*}
 \item[(b)] Since $\int_{\{|z|\le \vv\}} |z|^2\nu(\d z)<\infty$,
 \beg{equation*} \beg{split} & \int_{\{|z|\le \vv\}} \big|W(x+\si_2(x)z)-W(x)-\<\nn W(x), \si_2(x) z\>\big|\nu(\d z)\\
 &\le \ff 1 2 \sup_{|y|\le |x|+\|\si_2(x)\|} \|\nn^2 W(y)\| \int_{\{|z|\le \vv\}} |\si_2(x)z|^2 \nu(\d z)   \end{split}\end{equation*} is locally bounded in $x\in\R^d$. Noting that
 $$ \nn^2W(x)(v,v)  = \varphi'(|x|) \Big(\ff{|v|^2}{|x|}-\ff{\<x,v\>^2}{|x|^3}\Big)+\varphi''(|x|) \ff{\<x,v\>^2}{|x|^2}
 \le \ff{|v|^2(B(|x|)-|x|B'(|x|)}{B(|x|)^2(1+|x|)},$$ we obtain
  \beg{equation*} \beg{split} &  \int_{\{|z|\le \vv\}} \big[W(x+\si_2(x)z)-W(x)-\<\nn W(x), \si_2(x) z\>\big]\nu(\d z)\\
  &\le   \|\si_2(x)\|^2\tt B_\vv(x) \int_{\{|z|\le \vv\}} |z|^2\nu(\d z). \end{split}\end{equation*}
\item[(c)]  We have
   \beg{equation*}\beg{split}  & \inf_{\R^d} W- W(x)\le  W(x+\si_2(x)z)-W(x)\\
   &\le
     \varphi(|x+\si_2(x)z|)-\varphi(|x|)
     \le\int_{|x|}^{|x|+\|\si_2(x)\|\cdot |z|}\ff{\d s}{B(s)}. \end{split}\end{equation*}
 Since $\nu(\{|z|>\vv\})<\infty$ and $A_\vv<0$, this implies that $$\int_{\{|z|>\vv\}}  [W(x+\si_2(x)z)-W(x) ]\nu(\d z)$$ is locally bounded and
$$   \int_{\{|z|>\vv\}} \big[W(x+\si_2(x)z)-W(x)\big]\nu(\d z)
 \le \int_{\{|z|>\vv\}}\nu(\d z)\int_{|x|}^{|x|+\|\si_2(x)\|\cdot|z|}\ff{\d s}{B(s)}.$$   \end{enumerate}

By combining \eqref{AB} with (a)-(c), we conclude that $\scr LW$ is locally bounded satisfying \eqref{C3}.
\end{proof}

\beg{proof}[Proof of Corollary \ref{C3.2}] By Theorem \ref{T3.1}, for each situations it suffices to choose $B$ such that   one of
 \eqref{C1} and \eqref{C2} holds for some $\vv\in (0,\Theta^{-1})$.

{\bf Case (1)}.   We take $B(r)= (1+r)^\dd$ for some $\dd\in (1,\theta)\cap(0,1].$ Then $\int_0^\infty\ff{\d s}{B(s)} <\infty$ such that
$$\limsup_{|x|\to\infty}\bigg(\ff{\|\si_2(x)\|}{B(|x|)}+ \int_{\{|z|>\vv\}}\nu(\d z)\int_{|x|}^{\infty}\ff{\d s}{B(s)}\bigg)=0.$$ Next, since $\dd>1$, for any $\vv\in (0,\Theta^{-1})$ we have
$$\limsup_{|x|\to\infty} \tt B_\vv(x)\|\si_2(x)\|^2 \le  \limsup_{|x|\to\infty}\ff{\|\si_2(x)\|^2}{(1-\vv\Theta)^{1+\dd}|x|^{1+\dd}}=0.$$ Therefore, \eqref{C}   implies  \eqref{C1}.

{\bf Case (2)}.  We take $B(r)=1+r.$  Then $\int_0^\infty\ff{\d s}{B(s)} =\infty$  and by \eqref{20},
\beg{equation*}\beg{split}&\limsup_{|x|\to\infty}\bigg\{  \|\si_2(x)\|^2\tt B_\vv(x)\int_{\{|z|\le\vv\}} |z|^2\nu(\d z) +\ff{\|\si_2(x)\|\int_{\{\vv<|z|\le 1\}}|z|\nu(\d z)}{B(|x|)}\\
&\qquad \qquad \qquad\qquad +\int_{\{|z|>\vv\}} \nu(\d z)\int_{|x|}^{|x|+\|\si_2(x)\|\cdot |z|} \ff{\d s} {B(s)} \bigg\}  \\
& \le\ff{\vv^2\Theta^2\int_{\{|z|\le\vv\}}|z|^2\nu(\d z)}{2(1-\vv\Theta)^2} +\Theta \int_{\{\vv<|z|\le 1\}} |z|\nu(\d z)+\int_{\{|z|>\vv\}} \log (1+\vv\Theta |z|)\nu(\d z)<-D.\end{split}\end{equation*}  Thus, \eqref{C2} follows from \eqref{C}.

{\bf Case (3)}.    We take $B(r)= (1+ r)^\theta.$ Since $\|\si_2\|$ is bounded and $\theta\in (0,1)$, we have $\Theta=0$, $\int_0^\infty \ff{\d s}{B(s)}=\infty$ and
\beq\label{32}\lim_{|x|\to\infty}  \|\si_2(x)\|^2\tt B_\vv(x)=0.\end{equation} Moreover,
\beq\label{33}\int_{|x|}^{|x|+\|\si_2(x)\|\cdot|z|}\ff{\d s}{B(s)} \le \min\Big\{\ff{\|\si_2(x)\|\cdot |z|}{(1+|x|)^\theta},\ \ff{\|\si_2(x)\|^{1-\theta}|z|^{1-\theta}}{1-\theta}\Big\}.\end{equation}  Since $\|\si_2\|$ is bounded and $\int_{\{|z|\ge 1\}} |z|^{1-\theta}\nu(\d z)<\infty$, by \eqref{33} and the dominated convergence theorem we obtain
$$\limsup_{|x|\to\infty} \int_{\{|z|\ge 1\}} \nu(\d z)\int_{|x|}^{|x|+\|\si_2(x)\|\cdot |z|} \ff{\d s} {B(s)}=0.$$ Then \eqref{C2} with $\vv=1$ follows from \eqref{C}.

{\bf Case (4)}. We first observe that \eqref{21} implies
\beq\label{21'}\limsup_{|x|\to\infty}   \int_{\{|z|> 1\}}\ff{|z|\cdot \|\si_2(x)\|}{|x|^\theta\land (|x|+\|\si_2(x)\|\cdot |z|)^\theta}\,\nu(\d z)<-D.\end{equation} Since when $\theta\ge 0$ \eqref{21} is equivalent to \eqref{21'}, we only consider $\theta<0$. In this case, for any $s>1$ there exists a constant $C(s)>0$ such that
$$(|x|+\|\si_2(x)\|\cdot |z|)^{-\theta}\le s|x|^{-\theta}+C(s) (\|\si_2(x)\|\cdot |z|)^{-\theta}.$$ Moreover, by
 \eqref{21} and $\int_{\{|z|>1\}}|z|^{1-\theta}\nu(\d z)<\infty$, we have  $\|\si_2(x)\|^{1-\theta}\to 0$ as $|x|\to\infty$  so that
\beg{equation*}\beg{split} &\limsup_{|x|\to\infty}   \int_{\{|z|> 1\}}\ff{|z|\cdot \|\si_2(x)\|}{  (|x|+\|\si_2(x)\|\cdot |z|)^\theta}\,\nu(\d z)\\
&\le \limsup_{|x|\to\infty}   \bigg(\ff{s\|\si_2(x)\|}{|x|^\theta} \int_{\{|z|> 1\}}|z|\nu(\d z) + C(s)\|\si_2(x)\|^{1-\theta} \int_{\{|z|>1\}}|z|^{1-\theta}\nu(\d z)\bigg) \\
&=\limsup_{|x|\to\infty} \ff{s\|\si_2(x)\|}{|x|^\theta}  \int_{\{|z|> 1\}} |z| \nu(\d z),\ \ s>1.\end{split}\end{equation*} Since $s>1$ is arbitrary,
we conclude that \eqref{21} implies \eqref{21'}.

Now, we take $B(r)= (1+ r)^\theta.$  Then $\int_0^\infty \ff{\d s}{B(s)}=\infty, \tt B_\vv(x)={\rm O}((1+|x|)^{-(1+\theta)})$ for large $|x|$, and by   \eqref{21} we have $\Theta=0$. So, combining \eqref{21'} with \eqref{32} and \eqref{33}, we obtain
\beg{equation*}\beg{split}&\lim_{|x|\to\infty}\bigg\{  \|\si_2(x)\|^2\tt B_\vv(x) \int_{\{|z|\le 1\}}|z|^2\nu(\d z) +\int_{\{|z|> 1\}} \nu(\d z)\int_{|x|}^{|x|+\|\si_2(x)\|\cdot |z|} \ff{\d s} {B(s)} \bigg\}  \\
&\le \limsup_{|x|\to\infty} \int_{\{|z|> 1\}}\ff{|z|\cdot \|\si_2(x)\|}{|x|^\theta\land (|x|+\|\si_2(x)\|\cdot |z|)^\theta}\,\nu(\d z)<-D.\end{split}\end{equation*}Then  \eqref{C} implies \eqref{C2} with $\vv=1$.
\end{proof}

\section{Proofs of Theorems  \ref{T1.1} and \ref{T1.2}}

To apply Theorem \ref{T2.1}, we take   $b_t=b$ and $a_t(z)=\si z$ such that (\ref{E1}) reduces back to (\ref{E1.1}). In this case we have
$$\psi_s(z)= \si^{-1} (X_T^{s,x}(z)-X_T(x))$$ and for $s\in (0,T]$
$$\d \nn X_t^{s,x}= \{\nn b (X_t^{s,x})\} \nn X_t^{s,x}\,\d t +\si  \dd_s(\d t),\ \ \nn X_0^{s,x}=0.$$
Thus,  $\nn\psi_s(z)= \si^{-1} \nn X_T^{s,x}$ and for $t\ge s$,
$$\d(\si^{-1} \nn X_t^{s,x}) = \big\{\si^{-1} \nn b(X_t^{s,x}) \si\big\} \si^{-1} \nn X_t^{s,x}\d t,\ \ \si^{-1} \nn X_s^{s,x}= I.$$ Combining this with (\ref{A0}) we obtain
 \beq\label{D1} |{\rm det}
(\nn\psi_{s})^{-1}|=\ff 1 {|{\rm det}(\nn\psi_s)|}\le
\e^{-\ll_1(T-s)d},\end{equation} and
\beq\label{D2}\sup_{z\in\R^d\setminus\{0\}}
\ff{|z|}{|\psi_s^{-1}(z)|} = \sup_{z\in\R^d\setminus\{0\}}
\ff{|\psi_s(z)|}{|z|}\le \sup_{z\in\R^d\setminus\{0\}}
\ff{\int_0^1|\nn_z\psi_s(rz)|\d r }{|z|}\le  \e^{\ll_2
(T-s)}.\end{equation}

\

\beg{proof}[Proof of Theorem \ref{T1.1}] Let $\ff {\kk_1}{|z|^{d+\aa}}\le \ff{\nu(\d z)}{\d z}\le \ff {\kk_2}{|z|^{d+\aa}}$.
Then by \eqref{D1} and \eqref{D2} we obtain
\beg{equation*}\beg{split} & (\nu\circ\psi_s^{-1})(\d z)  \le \ff{\kk_2 |{\rm det}
\nn\psi_s^{-1}|(z)}{|\psi_s^{-1}(z)|^{d+\aa}}\,\d z
 =\ff{\kk_2 |{\rm det}(\nn\psi_s)^{-1}|(\psi_s^{-1}(z))}{|\psi_s^{-1}(z)|^{d+\aa}}\,\d z\\
&\le \ff{\kk_2 \e^{\ll_2(T-s)(d+\aa)-\ll_1(T-s)d}}{|z|^{d+\aa}}\,\d z
\le \ff{\kk_2}{\kk_1}\,\e^{\ll_2(T-s)(d+\aa)-\ll_1(T-s)d}
\nu(\d z).\end{split}\end{equation*} By Theorem \ref{T2.1}, this proves Theorem \ref{T1.1}(1).

Next,  to prove the existence of invariant probability measure using Corollary \ref{C3.2}, we take $ \si_1(x)=0$ and $\si_2(x)=\si.$ It is easy to see that $\Theta=0$ and $\int_{\{|z|\ge 1\}}\log (1+|z|)\nu(\d z)<\infty.$   Since    $\ll_2(d+\aa)-\ll_1d<0$  implies $\ll_2<0$, \eqref{C}  and \eqref{20} hold for $\theta=1.$ Then according to Corollary \ref{C3.2} for $\theta=1$,
  $P_t$ has an   invariant probability measure $\mu$. Moreover, since   $\ll_2<0$ implies
  $$\lim_{t\to\infty} |X_t(x)-X_t(y)|\le\lim_{t\to\infty} |x-y|\e^{\ll_2t}=0,\ \ x,y\in\R^d,$$
we conclude that  $P_tf\to\mu(f)$ as $t\to\infty$ holds for all $f\in C_b(\R^d).$ Thus, $\mu$ is the unique invariant probability measure of $P_t$.
Since $\Phi\in C^2((0,\infty))$, for $f\in C_b^2(\R^d)$ with $\inf f>0$ we have $\GG_{\Phi,\nu}(f)\in C_b(\R^d)$. By
letting $t\to\infty$ in the semigroup $\Phi$-entropy inequality in Theorem \ref{T1.1}(1), we prove (\ref{PE'}) for the desired constant $C$ and positive
$f\in C_b^2(\R^d)$ with $\inf f>0.$ By a simple
approximation argument, \eqref{PE'} holds for all $f\in \B_b^+(\R^d)$. 
\end{proof}

\beg{proof}[Proof of Theorem  \ref{T1.2}] (a) Let $\ll_2\le
0$ and $\rr$ be decreasing.    By (\ref{D2}) we have
$|z|\le |\psi_s^{-1}(z)|$, so that
\beq\label{AA}\rr(|\psi_s^{-1}(z)|)\le
\rr(|z|).\end{equation}Combining this with (\ref{D1}) and (\ref{D2})
we obtain
$$ (\nu\circ\psi_s^{-1})(\d z) \le \ff{\kk_2
|{\rm det}(\nn\psi_s)^{-1}|(z)\rr(|\psi_s^{-1}(z)|) }{|\psi_s^{-1}(z)|^{d+\aa}} \,\d z\le
\e^{\ll_2(T-s)(d+\aa)-\ll_1(T-s)d}\nu(\d z).$$ According to the proof of Theorem \ref{T1.1}, this proves the first  assertion in (I).

(b) By an approximation argument, for  \eqref{*W} we may assume that $f\in C_{c,+}^2(\R^d).$ We first consider the case that $b\in C^2(\R^d;\R^d)$ with bounded $\nn^2b$ and $\nu(1_{\{|\cdot|>1\}}|\cdot|)<\infty$. Then by the boundedness of  $\nn b$ (due to \eqref{A0}) and $\nn^2 b$, we see that $\|\nn X_t\|_\infty$ and $\|\nn^2X_t\|_\infty$ are locally bounded in $t\ge 0$, since for any $u,v\in\R^d$,
\beg{equation*}\beg{split} &\d \nn_u X_t= \nn b(X_t) \nn_u X_t\d t,\ \ \nn_u X_0=u,\\
&\d \nn_u\nn_v X_t= \big\{\nn^2 b(X_t) (\nn_u X_t,\nn_v X_t)+\nn b(X_t) \nn_u\nn_v X_t\big\}\d t,\ \ \nn_u\nn_v X_0=0.\end{split}\end{equation*}
This implies that   $P_t C_b^2(\R^d)\subset C_b^2(\R^d)$ for any $t\ge 0$ with $\|\nn P_tf\|_\infty$ and $\|\nn^2 P_tf\|_\infty$   locally bounded in $t.$ Next, by \eqref{A0} and $\nu(|\cdot|1_{\{|\cdot|\ge 1\}})<\infty$, it is easy to see that $W(x):= \ss{|x|^2+1}$ satisfies
\beq\label{W2} \scr LW(x)\le C_1-C_2 |x|,\ \ x\in\R^d\end{equation} for some constants $C_1,C_2>0.$ Thus, the invariant probability measure $\mu$ satisfies $\mu(|\cdot|)<\infty.$ Moreover, by the boundedness of $\nn b$, for any $f\in C_b^2(\R^d)$   there exists a constant $C_3>0$ such that $\scr Lf(x)\le C_3(1+|x|)$. So, by \eqref{W2} and $|\cdot|\le W$,
$$ \ff{|P_tf-f|}{t}\le \ff{C_3} t\int_0^t (1+\E |X_s|)\d s\le C_3 (1+W+C_1),\ \ t\in (0,1].$$Since the upper bound is integrable with respect to $\mu$, by the dominated convergence theorem  we obtain
$$\int_{\R^d} \scr L f\d\mu= \int_{\R^d}\lim_{t\to 0}\ff{P_tf-f}t\d\mu= \lim_{t\to\infty}\ff 1 t\int_{\R^d} (P_tf-f)\d\mu=0,\ \ f\in C_b^2(\R^d).$$
Since $P_t C_b^2(\R^d)\subset C_b^2(\R^d)$, for any $f\in C_{c,+}^2(\R^d)$ we have $\Phi(P_tf)\in C_b^2(\R^d)$, so that  $$\int_{\R^d} \scr L\Phi(P_tf)\d\mu=0.$$ Hence,  \eqref{W3} holds for $P_tf$ in place of $f$. Therefore, it follows from (\ref{PE'}) that
$$\ff{\d}{\d t}\Ent_\mu^\Phi(P_tf) =-\EE(\Phi'(P_tf),P_tf) =-\int_{\R^d} \GG_{\Phi,\nu}(P_tf)\d\mu \le -\ff 1 {C}\Ent_\mu^\Phi(P_tf),\ \ t\ge 0.$$ This implies \eqref{*W} for $f\in C_{c,+}^2(\R^d)$  since according to the first assertion \eqref{PE'} holds for  $C= \ff{\kk_2}{\kk_2(\ll_1 d-\ll_2(d+\aa))}.  $

(c) Assume that $\nu(1_{\{|\cdot|>1\}}|\cdot|)<\infty$. To apply the assertion proved in (b),  we make a standard regularization of $b$ as follows:
$$b_\vv(x)= \ff 1 {(\pi\vv)^{d/2}} \int_{\R^d} b(y)\e^{-|x-y|^2/\vv}\d y,\ \ x\in \R^d, \vv\in (0,1).$$
Since \eqref{A0} is equivalent to the dissipative property of $b(x)-\ll_2 x$ and $\ll_1 x-b(x)$, according to \cite[Theorem 9.19]{DR} we conclude that for every $\vv\in (0,1)$, $b_\vv\in C^2(\R^d;\R^d)$ with bounded $\nn^2 b_\vv$ and \eqref{A0} holds for $b_\vv$ in place of $b$. Then by (b), we have
\beq\label{W4} \Ent_{\mu_\vv}^\Phi(P^\vv_t f)\le \e^{-t/C} \Ent_{\mu_\vv}^\Phi(f),\ \ f\in \B_b^+(\R^d), t\ge 0,\end{equation} where $C= \ff{\kk_2}{\kk_2(\ll_1 d-\ll_2(d+\aa))},$ $P_t^\vv$ and $\mu_\vv$ are the semigroup and invariant probability measure for the equation
$$\d X_t^\vv= b(X_t^\vv)\d t +\si \d L_t,\ \ X_0^\vv=X_0.$$
Moreover, by the boundedness of $\nn b$,
$$|b_\vv(x)-b(x)|\le \ff{\|\nn b\|_\infty}{(\pi\vv)^{d/2}}\int_{\R^d}|x-y|\e^{-|x-y|^2/\vv}\d y\le c\ss\vv,\ \ x\in\R^d,\vv\in (0,1)$$ holds for some constant $c>0.$ Combining this with \eqref{A0} we obtain
\beg{equation*}\beg{split} \d |X_t-X_t^\vv|&=\ff{\<X_t-X_t^\vv,b(X_t)-b(X_t^\vv)\>+ \<X_t-X_t^\vv, b(X_t^\vv)-b_\vv(X_t^\vv)\>}{|X_t-X_t^\vv|}\d t\\
&\le \big\{\ll_2|X_t-X_t^\vv|+c\ss\vv\big\}\d t.\end{split}\end{equation*}Since $\ll_2<0$, this implies
$$|X_t-X_t^\vv|\le\ff{c\ss\vv}{-\ll_2}=:c'\ss\vv,\ \ t\ge 0,\vv\in (0,1).$$ Then for any $f\in C_b^1(\R^d)$,
\beq\label{W7} \|P_t^\vv f-P_tf\|_\infty \le\|\nn f\|_\infty c'\ss\vv,\ \ t\ge 0,\vv\in (0,1).\end{equation} Hence,
\beq\label{W8} |\mu_\vv(f)-\mu(f)|= \lim_{t\to\infty} |P_t^\vv f(0)- P_t f(0)|\le \|\nn f\|_\infty c'\ss\vv,\ \ f\in C_b^1(\R^d), \vv\in (0,1).\end{equation}
Combining  \eqref{W7} and \eqref{W8},   for any $f\in C_{b}^1(\R^d)$ with $\inf f>0$ we obtain
\beg{equation*}\beg{split} &\limsup\big|\mu_\vv(\Phi(P_t^\vv f))- \mu(\Phi(P_t f))\big|\\
&\le \limsup\big\{|\mu_\vv(\Phi(P_tf))-\mu(\Phi(P_tf))|
+\|\Phi(P_tf)-\Phi(P_t^\vv f)\|_\infty\big\}=0,\ \ t\ge 0.\end{split}\end{equation*}  Therefore,  letting $\vv\downarrow 0$ in \eqref{W4}, we prove \eqref{*W} for $f\in C_b^1(\R^d)$ with $\inf f>0,$ and thus also for $f\in \B_b^+(\R^d)$ by an approximation argument.

(d) In general, for any $n\ge 1$ let $\nu_n(\d z)= 1_{\{|z|\le n\}}\nu(\d z)$ and $\tt\nu_n= \nu-\nu_n.$ Write $L_t= L_t^n+\tt L_t^n$, where $L_t^n$ and $\tt L_t^n$ are L\'evy processes with L\'evy measures $\nu_n$ and $\tt\nu_n$ respectively. Consider the equation
$$\d X_t^n= b(X_t^n)\d t +\si \d L_t^n,\ \ X_0^n=X_0,$$ and let $P_t^n$ and $\mu_n$ be the associate semigroup and invariant probability measure respectively. Then by (c), we have
\beq\label{CD} \Ent_{\mu_n}^\Phi(P_t^nf)\le \e^{-t/C} \Ent_{\mu_n}^\Phi(f),\ \ t\ge 0, n\ge 1,f\in \B_b^+(\R^d).\end{equation} 
Noting that
$$\d (X_t-X_t^n)= (b(X_t)-b(X_t^n))\d t+\si\d\tt L_t^n,\ \ X_0-X_0^n=0,$$   by \eqref{A0} with $\ll_2<0$ we may find  a constant $\ll>0$ such that
such that
$$\d|X_t-X_t^n|^\vv\le \bigg(\int_{\{|z|>n\}}\big(|X_t-X_t^n+\si z|^\vv- |X_t-X_t^n|^\vv\big)\nu(\d z) -\ll|X_t-X_t^n|^\vv\bigg) \d t +M_t$$ holds for some local martingale $M_t$. Since 
$$\int_{\{|z|>n\}}\big(|X_t-X_t^n+\si z|^\vv- |X_t-X_t^n|^\vv\big)\nu(\d z)\le \int_{\{|z|>n\}}|\si z|^\vv\nu(\d z)=:\dd_n\downarrow 0\ \text{as}\ n\uparrow\infty,$$ we obtain
$$\E |X_t-X_t^n|^\vv\le \ff{\dd_n}\ll,\ \ t\ge 0.$$ Thus,
$$|P_t^nf-P_tf|\le (\|\nn f\|_\infty+2\|f\|_\infty)\E(1\land |X_t-X_t^n|)\le \ff{\dd_n}\ll (\|\nn f\|_\infty+2\|f\|_\infty),\ \ f\in C_b^1(\R^d).$$
Then, according to the argument in the end of (c) using this estimate to replace \eqref{W7}, we prove \eqref{*W} by letting $n\to\infty$ in \eqref{CD}.

(e) Let $\ll_1\ge 0,\theta_1\ge 1$ and $\rr$ be increasing. Then
$|z|\ge |\psi_s^{-1}(z)|$, so that (\ref{AA}) holds and the
remainder of the proof is similar to (a).
\end{proof}

\paragraph{\bf Acknowledgement.} The author would like to thank   Jian Wang for helpful comments.

\end{document}